\documentclass[a4 paper, 11pt]{article}

\usepackage[utf8]{inputenc}
\usepackage[T1]{fontenc}
\usepackage{amsmath, amsfonts, amssymb}
\usepackage{graphicx}
\usepackage{fancyhdr}
\usepackage[left=3cm, right=3cm, top=3cm, bottom=3cm]{geometry}
\usepackage{enumitem}
\usepackage{multirow}
\newtheorem{theorem}{Theorem}

\newtheorem{lemma}{Lemma}

\newenvironment{proof}{\par\noindent\textbf{Proof.}}{\hfill $\blacksquare$\par\smallskip}

\begin{document}

\begin{center}
\textbf{\LARGE{On a Family of Hypergeometric Polynomials}}\\
\end{center}

\begin{center}
\textrm{\normalsize{\textbf{Kikunga\  Kasenda\ Ivan}}}\\
\end{center}

\begin{center}
\textit{\scriptsize{Physics Department, Faculty of Sciences, University of Kinshasa,P.O.Box 190 Kin\ XI, Kinshasa, D.R.Congo.}}\\
\end{center}

\begin{center}
\textit{\scriptsize{Department of Engineering Sciences, Faculty of Sciences and Technologies, Université Loyola du Congo, P.O. Box 3724 Kin I, Kinshasa, D.R.Congo.}}\\
\end{center}

\begin{center}
{\scriptsize{Email : ivan.kikunga@unikin.ac.cd, ivan.kikunga@ulc-icam.com}}\\
\end{center}

\begin{abstract}
We work on the SCE problems. We establish the expressions of three integrals' sequences, related to it, in terms of five families of polynomials. Relations between these integrals are demonstrated and we focus on one of the three problems : the determination of the family of polynomials noted $e_n (n \in \mathbb{N})$. We show that these polynomials are hypergeometric. From this property, the NU method can be applied to this family. We have been able to determine the Rodrigues formula. These polynomials have properties that distinguish them from classical hypergeometric polynomials. We state and demonstrate the theorem adapted to the determination of the $e_n$ generating function. Finally, the sequence of polynomials studied is expressed in terms of associated Laguerre polynomials with negative upper indices.
\end{abstract}
\

\section{Introduction}
\

The problems of mathematical physics often lead to the resolution of mathematical equations of any form. The literature presents analytical resolution methods, approximate as well as numerical. This is where our problem lies. We will pose the problem of determining three sequences of indefinite integrals. These integral problems are related to the determination of some polynomial families. These are solutions of a differential equation of the second order of the form :
\begin{equation}
\label{eql01}
y''+P(x) y' + Q(x) y = R(x),
\end{equation}
where $P, Q, R$ are functions of the real variable $x$. In the case where these functions are polynomial or rational, a particular solution of \eqref{eql01} is a polynomial.\\
The Frobenius method makes it possible to obtain the families of polynomials sought, solutions of differential equations of the type \eqref{eql01}.\\
The detailed study of the problems posed reveals that these families of polynomials are linked and the study of only one of the cases is sufficient to determine the others.\\
The differential equation \eqref{eql01} has been the subject of many studies, for example H. Ciftci \textit{et al}\cite{Art03-14} and Y-Z. Zhang\cite{Art03-15} who were interested in the general conditions so that a homogeneous equation of the type \eqref{eql01} admits polynomial solutions as well as the underlying applications.\\

In some cases, the resolution of Eq.\eqref{eql01} leads to the study of special functions of Physics. This theory is an active branch of research in Theoretical and Mathematical Physics.\\
Literature often presents each special function as an isolated case and requires a detailed study on its own. These studies often lead to the establishment of Rodrigues formulas \cite{Art03-16}, recurrence and other relationships for each type of "special" polynomials (See \cite{Art03-05, Art03-09, Art03-10, Art03-11, Art03-12, Art03-13, Art03-07, Art03-08}).\\

The method of Nikiforov-Uvarov (NU Method)\cite{Art03-05, Art03-09, Art03-33} present all special functions as arising from a single theory : the theory of differential equations of the hypergeometric type\footnote{It is Eq.\eqref{eql01} homogeneous, with $P \ \text{and} \ Q$ polynomials or rational functions of some form.}. This method proposes a construction of special functions from a simple idea : to treat them from a single point of view. In some cases, the solutions of \eqref{eql01} are hypergeometric functions (See \cite{Art03-05, Art03-11, Art03-13, Art03-17, Art03-19, Art03-20}).\\
Note that the NU method has been used successfully to solve some theoretical and mathematical physics problems. For example :
\begin{itemize}
\item the resolution of the Schr\"{o}dinger's equation by B.I. Ita \textit{et al} \cite{Art03-21}-\cite{Art03-22}, A. D. Antia \textit{et al} \cite{Art03-23}, H. Louis \textit{et al} \cite{Art03-24}-\cite{Art03-25}, I. B. Okon \textit{et al} \cite{Art03-26}, M. Abu-Shady \cite{Art03-27}, C. A. Onate \textit{et al}\cite{Art03-02};

\item the resolution of Klein-Gordon's equation by B. I. Ita \textit{et al} \cite{Art03-28} and T. O. Magu \textit{et al} \cite{Art03-29};

\item some applications in Mathematical Physics by B. G\"{o}nul \textit{et al} \cite{Art03-31};

\item some applications in High Energy Physics, by exemple, M. Abu-Shady \textit{et al}\cite{Art03-32},...;

\item some applications of the "extended NU method" by H. Karayer \textit{et al} \cite{Art03-33} and C. Quesne \cite{Art03-34}.
\end{itemize}

In the present work, we generated a family of hypergeometric polynomials whose particular property will allow to present the adapted integral representation. The consequence is that, it leads to the correct expression of the generating function. Finally, we've found a link with the associated Laguerre polynomials with negative orders\footnote{We call \textit{order} of an associated Laguerre polynomials $L_n^{(\alpha)}$, 
the higher index $\alpha$ of these polynomials.}. Similar cases of studies of Laguerre polynomials of particular orders have been carried out by K. N. Boyadzieh \cite{Art03-35} and I. K. Khabibrakhmanov \textit{et al}\cite{Art03-39}.\\

The paper is structured as follows. In the second section, we begin with a brief overview of the NU method as presented in \cite{Art03-05} and \cite{Art03-09}.
The section 3 is devoted to the resolution of the problematic (integral problems). From these problems are deduced differential equations which, solved by the Frobenius method, lead to the explicit expressions of the families of polynomials. The section closes by determining the generating functions for each family of polynomials. In Section 4, we use the NU method to determine the Rodrigues formula and the generating function. The particular and adapted integral representation of the polynomials is discussed. Section 5 is devoted to the complete resolution of the second-order differential equation to which a family of polynomials determined in the second section obeys. From this study follows a relationship between the family of associated Laguerre polynomials with negative orders and that of the polynomials of the study. The last section deals with a generalization of the E problem. Finally, the paper ends with conclusions.

%%%%%%%%%%%%%%%%%%%%%%%%%%%%%%%%%%%%%%%%%%%%%%%%%%%%%%%%%%%%%%%%%%%%%%%%%%%%%%%%%%%%%%%%%%%%%%%%%%%%%

\section{Nikiforov-Uvarov Method or NU Method}
\

In the following lines, we will recall some essential elements of the NU method\cite{Art03-05} that will be useful in the study.\\

Consider the following second-order differential equation :
\begin{equation}
\label{eqa76}
u''+\frac{C(x)}{A(x)} \ u'+\frac{\tilde{C}(x)}{A^2(x)} \ u=0.
\end{equation}
where $u$ is a function of the real or complex variable $x$, $C$ is a polynomial of degree not greater than one, $A$ et $\tilde{C}$, are polynomials of degrees no greater than two.\\
Indeed, any differential equation of the second order having at most three singular points can be written in the form \eqref{eqa76} by using a suitable change of variable\cite{Art03-33}. \\
Consider the change of variable :
\begin{equation}
\label{eqa77}
u = \varphi (x) \ y,
\end{equation}
the starting equation takes the following form :
\begin{equation}
\label{eqa78}
y''+ \left( 2 \ \frac{\varphi'}{\varphi} + \frac{C}{A} \right)y' + \left( \frac{\varphi''}{\varphi} + \frac{C}{A} \ \frac{\varphi'}{\varphi} + \frac{\tilde{C}}{A^2} \right)y = 0.
\end{equation}
In order to simplify \eqref{eqa78}, the following substitutions are made :
\[\frac{\varphi'}{\varphi} =\frac{D}{A}, \ D = \frac{1}{2} \left( B - A \right),\]
\hspace{30mm}with $D$ and $B$ polynomials of degrees no greater than one.\\
The coefficient of $y'$ in \eqref{eqa78} becomes $\frac{B}{A}$ and the one of $y$ ; $\frac{F}{A^2}$ with
\[F = \tilde{C} + C \ D + D^2 + D' \ A - D \ A'.\]
The polynomials $B$ and $F$ have degrees, respectively, not greater than one and two. Previous results generate a class of transformations(See \eqref{eqa77}) leaving \eqref{eqa76} invariant. Indeed, the equation \eqref{eqa76} becomes :
\begin{equation}
\label{eqa83}
y''+\frac{B(x)}{A(x)} \ y'+\frac{F(x)}{A^2(x)} \ y=0.
\end{equation}
The choice of the function $ D $ is arbitrary, so we can choose its coefficients such as
\begin{equation}
\label{eqa84}
F(x)=\lambda \ A(x),
\end{equation}
\hspace{50mm}where $\lambda$ is a constant.\\
The equation \eqref{eqa83} takes, at this moment, the form below :
\begin{equation}
\label{eqa85}
A(x) \ y'' + B(x) \ y' + \lambda \ y=0.
\end{equation}
The equation \eqref{eqa85} is called \textit{second-order differential equation of the hypergeometric type} and its solutions are \textit{the hypergeometric functions}. If the solution of the equation\eqref{eqa85} is a polynomial, we will speak of a \textit{hypergeometric polynomial}.\\
The NU method proposes a procedure to turn a differential equation of the form \eqref{eqa76} to the form \eqref{eqa85} (See \cite{Art03-05}).\\

By restricting oneself to the case of hypergeometric polynomials, it is possible to determine the generalized Rodrigues formula for special functions. To do this, the NU theory begins with this important property : \textit{"the derivatives of any order of hypergeometric type functions are hypergeometric functions."} Thus, the derivative of the nth order of the Eq.\eqref{eqa85} is the following hypergeometric differential equation :
\begin{equation}
\label{eqzz01}
A(x) \ y_n'' + B_n(x) \ y_n' + \mu_n \ y_n=0,
\end{equation}
with $B_n=B+nA' ; \mu_n = \lambda + n B' + \frac{n\left(n-1\right)}{2}A''$ and $y_n(x)=\frac{d^n}{dx^n}y(x)$.\\
In the case of hypergeometric polynomials, we can have $\mu_n=0$, so (\cite{Art03-05},\cite{Art03-09}) :
\begin{equation}
\label{eqa106}
\lambda = \lambda_n \equiv -nB'-\frac{n\left(n-1\right)}{2}A''.
\end{equation}
Thus, the explicit expression of hypergeometric polynomials, of degree $n$, is :
\begin{equation}
\label{eqa117}
y(x) = \frac{\beta_n}{\rho(x)} \ \frac{d^n}{dx^n}\Big\lbrack A^n(x) \ \rho(x) \Big\rbrack, \ n=0, 1, 2, \cdots
\end{equation}
This is the \textit{Rodrigues formula} sought, with $\beta_n$ a normalization factor and the function $\rho$, a solution of the following differential equation :
\begin{equation}
\label{eqi01}
\left(A \rho\right)' = B \rho.
\end{equation}
The function $\rho$ allows to write Eq.\eqref{eqa85} in a self-adjoint form. For the differential equation \eqref{eqzz01}, the function $\rho_n$ allowing to write it in a self-adjoint form is such that :
\begin{equation}
\label{eqzz02}
\left(A \rho_n\right)' = B_n \rho_n.
\end{equation}
From Eq.\eqref{eqi01} and \eqref{eqzz02}, we have :
\begin{equation}
\label{eqzz03}
\rho_n(x) = A^n(x) \rho(x).
\end{equation}

The NU method provides a general expression of the generating functions for hypergeometric polynomials. The starting point is the integral representation of these polynomials in the complex plane. By replacing this representation in the definition of the generating functions, we obtain a relation allowing to find the generating function of any family of hypergeometric polynomials.\\

Let $\rho$ be a function solution of Eq.\eqref{eqi01} and let us replace $n$ by $\nu$ in \eqref{eqa106}. There is a function $\rho_\nu$ given by \eqref{eqzz03}, such that the function $u$ (solution of Eq.\eqref{eqa76}) given by : 
\begin{equation}
\label{eqb03}
u(x) = \int_{(C)} \frac{\rho_\nu(s)}{(s-x)^{\nu+1}} ds.
\end{equation}
That is true if :
\begin{enumerate}
\item by calculating the derivatives $u'$ and $u''$, derivation and integration operations can be swapped :
\begin{eqnarray}
\begin{split}
\label{eqb08}
u'(x) = \left(\nu+1\right) \int_{(C)} \frac{\rho_\nu(s)}{(s-x)^{\nu+2}} ds,\\
u''(x) = \left(\nu+1\right)\left(\nu+2\right) \int_{(C)} \frac{\rho_\nu(s)}{(s-x)^{\nu+3}} ds,
\end{split}
\end{eqnarray}
\item the integration contour $(C)$ is chosen such that
\begin{equation}
\label{eqb09}
\frac{A^{\nu+1}(s)\rho(s)}{\left(s-x\right)^{\nu+2}}\bigg\vert^{s_2}_{s_1}=0,
\end{equation}
where $s_1$ and $s_2$ are the ends of the $(C)$ contour.
\end{enumerate}
Then, the equation \eqref{eqa85} admits particular solutions of the form :
\begin{eqnarray}
\begin{split}
\label{eqb04}
y_\nu(x) &= \frac{\gamma_\nu}{\rho(x)} \ \int_{(C)} \frac{A^\nu(s)\rho(s)}{(s-x)^{\nu+1}} ds,\\
&\gamma_\nu=\frac{\nu!}{2\pi i}\beta_\nu.
\end{split}
\end{eqnarray}
For $\nu$ an integer, $y_\nu$ is a polynomial. When $\nu$ is not an integer, Eq.\eqref{eqb04} still gives a particular solution of Eq.\eqref{eqa85}.\\
Let $F(x, t)$ be the generating function of the polynomials given by \eqref{eqb04}, we have :
\begin{equation}
\label{eqb06}
F(x,t)=\sum_{n=0}^{+\infty} \frac{1}{n!} \tilde{y}_n (x) t^n,
\end{equation}
where the polynomials $\tilde{y}_n$ are obtained by putting in \eqref{eqa117}, $\beta_n=1$.\\
Finally, the general expression of the generating functions of hypergeometric polynomials, given in \cite{Art03-05}, is 
\begin{equation}
\label{eqb07}
F(x,t)=\frac{\rho(s)}{\rho(x)} \frac{1}{1-A'(s)t}\Big\vert_{s=\xi(x,t)},
\end{equation}
with $s=\xi(x,t)$, solution of the equation $s-x-A(s)t=0$.

%%%%%%%%%%%%%%%%%%%%%%%%%%%%%%%%%%%%%%%%%%%%%%%%%%%%%%%%%%%%%%%%%%%%%%%%%%%%%%%%%%%%%%%%%%%%%%%%%%%%%

\section{Generation of Polynomials' Families Starting from a General Problematic}

\subsection{The General Problematic}

Let us consider the following integrals :

\begin{equation}
\label{eqa04}
S_n(x)=\int x^n \sin x \ dx,
\end{equation}
\begin{equation}
\label{eqa05}
C_n(x)=\int x^n \cos x \ dx, 
\end{equation}
\begin{equation}
\label{eqa06}
E_n(x)=\int x^n e^x \ dx, 
\end{equation}
with $n$, a natural integer. 
In the following, we will respectively designate these problems by the acronym SCE. There are tables showing the SCE expressions, see for example \cite{Art03-07} and \cite{Art03-08}. In most documents presenting these results, recurrence relations are often mentioned, which by iteration make it possible to generate the sequence of solutions. This method becomes tedious when the exponent $n$ becomes very large. Hence, the importance of having general results to determine SCE integrals for any value of $ n $, no matter how large. We present the results of the SCE integrals in the form of a theorem and we give the proof.

\begin{theorem}
Let the sequences above : \eqref{eqa04}, \eqref{eqa05} and \eqref{eqa06}. We can show that they can, respectively, take the following forms:
\begin{equation}
\label{eqa07}
S_n(x)=s_n(x) \ \cos{x} + \hat{s}_{n-1}(x) \sin{x} + \alpha_n,
\end{equation}
\begin{equation}
\label{eqa08}
C_n(x)=c_n(x) \ \sin{x} + \hat{c}_{n-1}(x) \ \cos{x} + \beta_n, 
\end{equation}
\begin{equation}
\label{eqa09}
E_n(x)=e_n(x) \ e^x + \gamma_n, 
\end{equation}
where $c_n, s_n \ \text{and} \ e_n$ are polynomials of degrees $(n \in \mathbb{N}); \hat{s}_{n-1} \ \text{and} \ \hat{c}_{n-1}$ polynomials of degrees $n-1$ with $n \in \mathbb{N^\star}$.\\
The $\alpha_n,\beta_n$ and $\gamma_n$ are, respectively, integration constants linked to each problem of order $n$.\\
By hypothesis, we must have $s_{-n} = \hat{s}_{-n} = c_{-n} = \hat{c}_{-n} = e_{-n}=0$, $\forall n \in \mathbb{N^{\star}}.$
\end{theorem}
The polynomials $s, \hat{s}, c, \hat{c} \ \text{and} \ e$ are to be determined.\\

\proof{To this end, we will use inductive proof on the "S problem", the other proofs can be deduced from it by an adequate change of variables $x = \frac{\pi}{2} \pm y$ or $x = iy$, respectively for "C" and "E problems". Apart from the variable change, the S and C problems can be related. It is enough to carry out an integration by parts, one obtains :
\begin{equation}
\label{eqa45}
S_n(x)=-x^n \cos x + n C_{n-1}(x),
\end{equation}
\begin{equation}
\label{eqa48}
C_n(x)=x^n \sin x - n S_{n-1}(x).
\end{equation}

Let's show that the proposal \eqref{eqa07} is true for $n=0$. An integration of \eqref{eqa04}, for $n=0$, gives $S_0(x) = - \cos x + cst.$ This corresponds to \eqref{eqa07} with $s_0(x)=-1$, $\hat{s}_{-1}(x)=0$ and $\alpha_0 = cst$. The polynomial $s_0$ is of degree zero and the associated integration constant is of zero order.\\
Suppose that the proposition \eqref{eqa07} is true for the index $m$. Let's try to establish its veracity for $m + 1.$ From \eqref{eqa07}, we have to check that
\begin{equation}
\label{eqb10}
S_{m+1}(x) = s_{m+1}(x) \cos x + \hat{s}_m (x) \sin x + \alpha_{m+1}.
\end{equation}
The polynomials $s_{m+1}$ and $\hat{s}_m$, being of degrees respectively $m+1$ and $m$.\\
From \eqref{eqa04} and \eqref{eqa05}, we have
\begin{equation}
\label{eqb11}
dS_m = x^m \ \sin x \ dx \  \text{and} \ dC_m = x^m \ \cos x \ dx.
\end{equation}
We can write
\[dS_{m+1} = x \ dS_m = x\left(s'_m + \hat{s}_{m-1}\right) \cos x \ dx + x\left(-s_m + \hat{s}'_{m-1}\right) \sin x \ dx.\]
Let $\lambda_m = x\left(s'_m + \hat{s}_{m-1}\right)$ and $\mu_{m+1}=x\left(-s_m + \hat{s}'_{m-1}\right)$. The polynomials $\lambda_m \ \text{and} \ \mu_{m+1}$ are of degrees respectively $m$ and $m+1$. We then write
\begin{equation}
\label{eqb13}
dS_{m+1} = \sum_{k=0}^m \tilde{\lambda}_m(k) \ dC_k  + \sum_{k=0}^{m+1} \tilde{\mu}_{m+1}(k) \ dS_k,
\end{equation}
with $dS_k$ and $dC_k$ given by \eqref{eqb11}.\\
The real $\tilde{\lambda}_m(k)$ are the order-$k$ coefficients of the polynomial $\lambda_m$ and the $\tilde{\mu}_{m+1}(k)$ are the order-$k$ coefficients of the polynomial $\mu_{m+1}.$\\
By integrating \eqref{eqb13}, we have:
\begin{equation}
\label{eqb13b}
S_{m+1}(x) = \sum_{k=0}^m \tilde{\lambda}_m(k) \ C_k(x)  + \sum_{k=0}^{m+1} \tilde{\mu}_{m+1}(k) \ S_k(x) + K,
\end{equation}
where $K$ is an integration constant.\\
In order to express the relation \eqref{eqb13b} as a function of $S_n$, let us use the relation \eqref{eqa48}, which finally allows us to write

\begin{equation}
\label{eqb14}
S_{m+1}(x) = \tau_{m+1}(x) \cos x + \kappa_m(x) \sin x + \theta_{m+1}. 
\end{equation}
with
\[\tau_{m+1}(x) = -\sum_{k=0}^m k \ \tilde{\lambda}_m(k) \ s_{k-1}(x) + \sum_{k=0}^{m+1} \tilde{\mu}_{m+1}(k) \ s_k (x),\]
\[\kappa_m(x) = \sum_{k=0}^m \tilde{\lambda}_m(k) \ \Big\lbrack x^k - k \ \hat{s}_{k-2}(x) \Big\rbrack + \sum_{k=0}^{m+1} \tilde{\mu}_{m+1}(k) \ \hat{s}_{k-1}(x)\]
and the integration constant of $m + 1$ order taken as
\[\theta_{m+1}=-\sum_{k=0}^m k \ \tilde{\lambda}_m(k) \ \alpha_{k-1} + \sum_{k=0}^{m+1} \tilde{\mu}_{m+1}(k) \ \alpha_k + K.\]
The polynomials $\tau_{m+1}$ and $\kappa_m$, are respectively of degree $m + 1$ and $m$. Of course $\tau_{m+1}=s_{m+1}, \kappa_m=\hat{s}_m$ and $\theta_{m+1}=\alpha_{m+1}$ to recover Eq.\eqref{eqb10}. This proves theorem 1.}

\subsection{Differential equations related to "SCE problems"}
\

The SCE problems can be put into differential form, in the sense that one can substitute for each of them a differential equation. The problem can be taken in another sense, that of the resolution of a differential equation. We will explicitly make this transformation for the S problem.\\

By deriving the two members of \eqref{eqa04} and \eqref{eqa07} and identifying them, we obtain the following relations :
\[s'_n(x) + \hat{s}_{n-1}(x)=0 \ \text{and} \ \hat{s}'_{n-1}(x) - s_n(x) = x^n,\]
then by eliminating from these relations the polynomials $\hat{s}_{n-1}$, we obtain the differential equation to which the polynomials $s_n$ must obey :
\begin{equation}
\label{eqa01}
s''_n+s_n=-x^n.
\end{equation}
Thus, the S problem is reduced to the determination of a single family of polynomials $\lbrace s_n; n \in \mathbb {N} \rbrace $. The same procedure applied to C and E problems leads to the following differential equations :
\begin{equation}
\label{eqa02}
c''_n + c_n = x^n, 
\end{equation}
and
\begin{equation}
\label{eqa03}
e'_n+e_n=x^n. 
\end{equation}
The polynomials $\hat{s}_{n-1}$ and $\hat{c}_{n-1}$ are, in fact, respectively bound to the polynomials $s_n$ and $c_n$. Thus, the number of unknown polynomials drops from five $\left(s_n, \hat{s}_{n-1}, c_n, \hat{c}_{n-1}, e_n\right)$ to three $\left(s_n, c_n, e_n\right)$. The $\hat{s}_{n-1}$ and $\hat{c}_{n-1}$ can be deduced from the following relationships : 
\begin{eqnarray}
\begin{split}
\label{eqa10}
\hat{s}_{n-1}(x)&=-s'_n(x),  \\
\hat{c}_{n-1}(x)&=c'_n(x).  \\
\end{split}
\end{eqnarray}
We can apply the same procedure to the equation \eqref{eqa02} but that is useless because the polynomials $s_n$ and $c_n$ are linked. To see it, we can replace \eqref{eqa07} in \eqref{eqa48} then by virtue of \eqref{eqa08}, we obtain :
\begin{eqnarray}
\begin{split}
\label{eqa33}
c_n &= x^n + n s'_{n-1},\\
c'_n &= -n s_{n-1}.
\end{split}
\end{eqnarray}
Other relationships can be found between the families $ \lbrace s_n; n \in \mathbb {N} \rbrace$ and $\lbrace c_n; n \in \mathbb {N} \rbrace$. The first group of these relations is obtained by using the relation of recurrence \eqref{eqa45} and by substituting the relations \eqref{eqa07} and \eqref{eqa08}, one finds:
\begin{eqnarray}
\begin{split}
\label{eqa46}
s_n &=-x^n + n \ \hat{c}_{n-2},\\
\hat{s}_n &=\left(n+1\right)c_n.
\end{split}
\end{eqnarray}
The second group is obtained by integrating once again \eqref{eqa45} by parts, this leads to
\begin{equation}
\label{eqc01}
S_n(x)=-x^n \ \cos x + n x^{n-1} \ \sin x -n\left(n-1\right)S_{n-2}(x),
\end{equation}
then substituting \eqref{eqa07}, we have the second group 
\begin{eqnarray}
\begin{split}
\label{eqa47}
s_n &= - x^n - n\left(n-1\right) s_{n-2},\\
\hat{s}_n &=\left(n+1\right)x^n-n\left(n+1\right)\hat{s}_{n-2}.
\end{split}
\end{eqnarray}
The third group is obtained by the same procedure, one uses the relation \eqref{eqa48} and the relations \eqref{eqa07} and \eqref{eqa08}. After identification, we find
\begin{eqnarray}
\begin{split}
\label{eqa49}
c_n &=x^n - n \hat{s}_{n-2},\\
\hat{c}_n &=-\left(n+1\right)s_n.
\end{split}
\end{eqnarray}
The last group is obtained by the same procedure as the second. The second partial integration on \eqref{eqa48} gives
\begin{equation}
\label{eqc02}
C_n(x)=x^n \ \sin x + n x^{n-1} \ \cos x -n\left(n-1\right)C_{n-2}(x).
\end{equation}
By virtue of \eqref{eqa07}, \eqref{eqa08} and by identification, we have
\begin{eqnarray}
\begin{split}
\label{eqa50}
c_n&=x^n - n\left(n-1\right) c_{n-2},\\
\hat{c}_n&=\left(n+1\right)x^n-n\left(n+1\right)\hat{c}_{n-2}.
\end{split}
\end{eqnarray}
Thus, taking into account all the recurrence relations obtained, it is possible to determine any polynomial of the S and C problems and to deduce the others. We have specifically chosen to treat problem E separately. This one can be conceived as a linear combination of the two others, via the relation of Euler. Thus, it is possible to reduce the determination of the polynomials $s_n$ (and the others related to it) to the sole search of the polynomials $e_n$.

\subsection{Explicit Determination of the $e_n$ Polynomials}
\

The $e_n$ are polynomial solutions of the differential equation \eqref{eqa03}; a particular solution. Let's find the general solution. We begin by writing Eq.\eqref{eqa03} in this form :
\begin{equation}
\label{eqz01}
y'+y=x^n,
\end{equation}
where $y$ is a function of the variable $x$ and $n \in \mathbb{N}$.\\
The homogeneous solution $y_h$ of Eq.\eqref{eqz01} is 
\begin{equation}
\label{eqz02}
y_h = K \ e^{-x},
\end{equation}
where $K$ is a integration constant.\\
To find the general solution of Eq.\eqref{eqz01}, we use the method of variation of the integration constant. The solution \eqref{eqz02} is then replaced in \eqref{eqz01} and we obtain: 
\begin{equation}
\label{eqz03}
\frac{d}{dx}K(x)=x^n \ e^{x}.
\end{equation}
We note that the determination of the function $K$ is the E problem (See \eqref{eqa06},\eqref{eqa09}). So, the general solution of Eq.\eqref{eqz01} is :
\begin{equation}
\label{eqz06}
y(x) = e_n(x) + C \ e^{-x},
\end{equation}
with $C$, an integration constant.\\

After this, let's find the polynomial solutions of Eq.\eqref{eqa03}. In order to find them, let's use the adapted Frobenius method, writing :
\begin{equation}
\label{eqa14}
e_n(x)=\sum_{l=0}^{n} a_l \ x^l,
\end{equation}
since we know $e_n$ are polynomials of degree $n$(See Theorem 1).\\
Replacing \eqref{eqa14} in \eqref{eqa03}, we have :
\begin{equation}
\label{eqa14b}
a_n \ x^n + \sum_{l=0}^{n-1} \big\lbrack \left(l+1\right)a_{l+1} + a_l \big\rbrack x^l = x^n.
\end{equation}
The expression \eqref{eqa14b} shows, as one would expect, that there exists an arbitrary constant $a_0$ and that the coefficient $a_n=1$. In addition, one obtains the following recurrence formula
\begin{equation}
\label{eqa60}
a_l=\frac{-a_{l-1}}{l} \ \text{for} \ l \in \lbrack 1 , n \rbrack,
\end{equation}
which can take the form
\begin{equation}
\label{eqa61}
a_l=\frac{\left(-1\right)^l}{l!}a_0.
\end{equation}
Since $a_n = 1$, we can find the independent term of the polynomials $e_n$
\begin{equation}
\label{eqa62}
a_0=n!(-1)^n.
\end{equation}
Thus, the coefficient of the term in $x^l$ is
\begin{equation}
\label{eqa63}
a_l=\left(-1\right)^{n+l}\frac{n!}{l!}.
\end{equation}
Finally, the $e_n$ have the form :
\begin{equation}
\label{eqa64}
e_n(x)=x^n+\sum_{l=0}^{n-1}\left(-1\right)^{l+n}\frac{n!}{l!}x^l.
\end{equation}
The polynomials $e_n$ obey the following recursion relations :
\begin{eqnarray}
\begin{split}
\label{eqe01}
e_n(x)&=x^n-n \ e_{n-1}(x);\\
e'_n(x)&= n \ e_{n-1}(x).
\end{split}
\end{eqnarray}
Since the polynomials $e_n$ are completely determined, they serve to deduce the polynomials of the other two problems. To do this, in \eqref{eqa06}, we successively put $x = iy$ and $x = -iy$. The obtained expressions are linear combinations of the $S_n$ and $C_n$, this allows, by Theorem 1 to express the following relations :
\begin{eqnarray}
\label{eqa68}
c_n (x) &= \frac{i^n}{2} \ \Big\lbrack \left(-1\right)^n e_n(ix) + e_n (-ix) \Big\rbrack,\\ \nonumber
c'_n (x) &= \frac{i^{n+1}}{2} \ \Big\lbrack \left(-1\right)^{n+1} e_n(ix) + e_n (-ix) \Big\rbrack,
\end{eqnarray}
and
\begin{eqnarray}
\label{eqa69}
s_n (x) &= \frac{i^n}{2} \ \Big\lbrack \left(-1\right)^{n+1} e_n(ix) - e_n (-ix) \Big\rbrack,\\ \nonumber
s'_n (x) &= \frac{i^{n+1}}{2} \ \Big\lbrack \left(-1\right)^n e_n(ix) - e_n (-ix) \Big\rbrack.
\end{eqnarray}
Moreover, by adding the first relations of \eqref{eqa68} and \eqref{eqa69} and their second relations, we see that
\begin{eqnarray}
\label{eqe02}
c_n(x)&=-s_n(x); \ \forall n \in \mathbb{N},\\
\hat{s}_n(x)&=\hat{c}_n(x); \ \forall n \in \mathbb{N}.
\end{eqnarray}
The relations \eqref{eqa69} make it possible to determine the polynomials $s_n$ and it is found that the results obtained are in agreement with their explicit expressions obtained by the Frobenius method. We have :
\begin{eqnarray}
\begin{split}
\label{eqa23}
s^{even}_{n}(x)&=-x^n+\sum_{l=0}^{n-2}\left(-1\right)^{\frac{l+n}{2}+1}\frac{n!}{l!}x^l, \ n= 2, 4, \cdots \\
s^{odd}_{n}(x)&=-x^n+\sum_{l=0}^{n-2}\left(-1\right)^{\frac{l+n}{2}}\frac{n!}{l!}x^l, \ n= 3, 5, \cdots \\
\end{split}
\end{eqnarray}
with \[s_0(x)=-1 \ \text{and} \ s_1(x)=-x.\]
At this point, we report that the SCE problems were treated by G. Dattoli \textit {et al} \cite{Art03-36} by the method of the negative derivation operator. The polynomials $e_n$ can be compared with the results of G. Dattoli \textit{et al} in \cite{Art03-37, Art03-38} by their study of the polynomials obtained by the truncated-exponential. We can also mention, for example, the study of the truncated-exponential by H. M. Srivastava \textit{et al} in \cite{Art03-46}, ...

\subsection{The Generating functions}
\

Let us denote the generating functions of the polynomials $s_n$, $c_n$ and $e_n$, respectively by $S$, $C$ and $E$. We will start by determining the generating function $S$. The generating function $C$ can be deduced from $S$, by virtue of the relation \eqref{eqe02}. Next, we will determine the generating function $E$. Finally, the relations between these polynomials make it possible to link the three generating functions.\\

The generating function of polynomials $s_n$, denoted $S$, is defined by :
\begin{equation}
\label{eqe03}
S(x,t)=\sum_{n=0}^{+\infty} \frac{1}{n!} \ s_n(x) \ t^n.
\end{equation}
To find its explicit form, we start from the differential equation \eqref{eqa01}. This multiplied by $\frac{1}{n!} \ t^n$, allows to obtain the differential equation to which the generating function must obey : 
\begin{equation}
\label{eqa53}
\frac{\partial^2 S}{\partial x^2} + S = -e^{xt}.
\end{equation}
The general solution of \eqref{eqa53} is 
\begin{equation}
\label{eqa54b}
S(x,t)=f(t) \cos x + g(t) \sin x - \frac{e^{xt}}{1+t^2},
\end{equation}
where $f$ and $g$ are arbitrary functions. For a polynomial solution of Eq.\eqref{eqa01}, the generating function is the particular solution of Eq.\eqref{eqa53} :
\begin{equation}
\label{eqa54}
S(x,t)=\frac{-e^{xt}}{1+t^2}.
\end{equation}
The generating function of the $c_n$ polynomials, denoted $C$, is defined by :
\begin{equation}
\label{eqa55}
C(x,t)=\sum_{n=0}^{+\infty}\frac{1}{n!} \ c_n(x) \ t^n,
\end{equation}
can be determined by the same procedure applied to the differential equation \eqref{eqa02}. A faster method is to use the relation \eqref{eqe02}. This implies that $S(x,t)=-C(x,t).$ Then 
\begin{equation}
\label{eqa57}
C(x,t)=\frac{e^{xt}}{1+t^2}.
\end{equation}
Let us now determine the generating function of $e_n$. Indeed,
\begin{equation}
\label{eqa70}
E(x,t)=\sum_{n=0}^{+\infty} \frac{1}{n!} e_n(x) t^n.
\end{equation}
The differential equation \eqref{eqa03} gives the equation
\begin{equation}
\label{eqa71}
\frac{\partial E}{\partial x} + E = e^{xt},
\end{equation}
whose a particular solution is
\begin{equation}
\label{eqa72}
E(x,t)=\frac{e^{xt}}{1+t}.
\end{equation}
The generating functions are connected by the relation
\begin{equation}
\label{eqa73}
-2S(x,t) = 2 C(x,t)=E(ix,-it)+E(-ix,it).
\end{equation}

%%%%%%%%%%%%%%%%%%%%%%%%%%%%%%%%%%%%%%%%%%%%%%%%%%%%%%%%%%%%%%%%%%%%%%%%%%%%%%%%%%%%%%%%%%%%%%%%%%%%%%%%

\section{Application of the NU Method to the $e_n$ Polynomials}
\

The polynomials $e_n$ are among the solutions of the differential equation \eqref{eqa03} which is a first-order equation. In this section, we will show that these same polynomials are also solutions of a second-order differential equation of the hypergeometric type and we will study the consequences of this property.

\subsection{The $e_n$ Polynomials Are Hypergeometric}
\

Consider the differential equation \eqref{eqa03}, which admit $e_n$ as particular solutions. By differentiating it, we obtain the equation $e''_ n + e'_n = nx^{n-1}.$ By multiplying the new equation of the second order by the variable $x$ and then subtracting the equation \eqref{eqa03} multiplied by the integer $n$, we obtain the following differential equation:
\begin{equation}
\label{eqa75}
x e''_n + \left(x-n\right) e'_n - n e_n =0.
\end{equation}
This equation is of hypergeometric type (See \cite{Art03-05}), it is of the same type as \eqref{eqa85} with $A(x) = x, B(x)= x-n$ and $\lambda = -n $. This allows us to use the NU method in the study of this family of polynomials.

\subsection{Rodrigues Formula for the $e_n$ Polynomials}
\

The NU method provides the general Rodrigues formula suitable for all hypergeometric polynomials (see Eq.\eqref{eqa117}). The $e_n$ being hypergeometric, this formula will be used. This passes by the determination of the function $\rho$, solution of the differential equation \eqref{eqi01}. For the family of polynomials $e_n$, the function $\rho$ is
\begin{equation}
\label{eqa120}
\rho(x)=x^{-n-1} \ e^x.
\end{equation}
The relation \eqref{eqa120} replaced in \eqref{eqa117} gives the Rodrigues formula sought :
\[e_n (x) = \beta_n \ x^{n + 1} \ e^{-x} \ \frac{d^n}{dx^n} \left(x^{-1} \ e^x \right)\]
with the constants $\beta_n$ to be determined. To do this, we use the fact that
$e_n(0)= a_0 = n! \ \left(-1\right)^n.$ We have to determine the constant term provided by \eqref{eqa117} and compare it to $e_n(0)$. The binomial formula applied to the nth derivative in \eqref{eqa117} makes it possible to rewrite this expression in the form :

\begin{equation}
\label{eqe05}
e_n(x) = \beta_n \ x^n \sum_{k=0}^n \ \left(-1\right)^k \ \frac{n!}{\left(n-k\right)!} \ x^{-k},
\end{equation}
the constant term is obtained for $k=n$, in \eqref{eqe05}. So,
\begin{equation}
\label{eqa122}
e_n(0) = n! \left(-1\right)^n \beta_n \Longrightarrow \beta_n=1.
\end{equation}
Finally, the Rodrigues formula for the $e_n$ is

\begin{equation}
\label{eqa123}
e_n(x) = x^{n+1} \ e^{-x} \ \frac{d^n}{dx^n}\left(x^{-1} \ e^x \right).
\end{equation}
One can, at this point, find an expression of the polynomials $s_n$ by using the relation of Rodrigues \eqref{eqa123}, one obtains :
\begin{equation}
\label{eqa125}
s_n(x) = \frac{-i^n}{2} \ x^{n+1} \Big\lbrack \left(-1\right)^n \ e^{-ix} \ \frac{d^n}{dx^n}\left(x^{-1} \ e^{ix} \right) + e^{ix} \ \frac{d^n}{dx^n}\left(x^{-1} \ e^{-ix} \right) \Big\rbrack.
\end{equation}
The relation \eqref{eqa125} was found by combining \eqref{eqa123} and \eqref{eqa69}.

\subsection{Generating Function Of the $e_n$}
\

The integral representation \eqref{eqb04} makes it possible to determine the generating functions of the hypergeometric polynomials. Indeed, the NU method proposes the relation \eqref{eqb07}. The case of the $e_n$ polynomials is a particular one. The determination of the generating function, by the NU method, requires sustained attention. Then, although hypergeometric, the $e_n$ polynomials do not have the same behavior as the classical hypergeometric polynomials known in the literature. Indeed, by observing the proof of \eqref{eqb07}, provided by \cite{Art03-05}, the numerator of \eqref{eqb04}, contains the expression \eqref{eqzz03}, such as the definition \eqref{eqb06} used provides the result \eqref{eqb07}. The functions $\rho_n$, in the case of the polynomials $e_n$, are independent of the integer $n$. The relation \eqref{eqb07} can not be applied because the conditions of its use are not fulfilled. Thus, it is necessary to find an expression exactly reproducing the generating function of the $e_n$, taking into account this peculiarity : the case where the function $\rho_n$ would be independent of the integer $n$. In the case of the classical hypergeometric polynomials, the problem does not arise. Thus, it is necessary to write a result, similar to that of Nikiforov-Uvarov and adapted to hypergeometric polynomials having the same peculiarity as the $e_n$. The fact that the $\rho_n$ function does not depend on $n$, leads us to pose that $\rho_n(x)=\sigma(x)$.\\

Before stating the new result, let us study the behavior of the function $\sigma$, in the case of differential equations of the hypergeometric type of the same kind as \eqref{eqa75}. Let us consider the hypergeometric differential equation \eqref{eqa85}, with
\begin{equation}
\label{eqf01}
A(x)=\alpha x+\beta ; \ B(x)=\gamma x +\delta.
\end{equation}
The $\rho$ function, in this case is given by :
\begin{equation}
\label{eqf02}
\rho(x)=\left(\alpha x+ \beta\right)^{\frac{-\beta \gamma + \alpha \left(\delta - \alpha\right)}{\alpha^2}} \ \exp{\left(\frac{\gamma}{\alpha}x\right)},
\end{equation}
and the function $\sigma$ by
\begin{equation}
\label{eqf03}
\sigma(x)=\left(\alpha x+ \beta\right)^{n+\frac{-\beta \gamma + \alpha \left(\delta - \alpha\right)}{\alpha^2}} \ \exp{\left(\frac{\gamma}{\alpha}x\right)}.
\end{equation}
The condition for not using the NU formula \eqref{eqb07}, in cases similar to the $e_n$, is
\begin{equation}
\label{eqf04}
n+\frac{-\beta \gamma + \alpha \left(\delta - \alpha\right)}{\alpha^2}= constant.
\end{equation}
Applied to Laguerre polynomials(or associated Laguerre) for
\[\alpha = 1, \beta = 0, \gamma = -1 \ \text{and} \ \delta = n \ \left(\text{or} \ m + 1 \right),\]
the condition \eqref{eqf04} is not verified, therefore the generating function can be determined by \eqref{eqb07}. In the case of $e_n$, it is not difficult to check that \eqref{eqf04} is respected. This is also the case, for example, for $\alpha=1, \beta=1, \gamma=-1 \ \text{and} \ \delta=-n$. This is the case of the family of polynomials $\lbrace \left(-1\right)^n \ \left(x + 1\right)^n ; n \in \mathbb{N} \rbrace$.\\

With all these previous considerations, we have to rewrite the integral representation of the $e_n$ (See Eq.\eqref{eqb03} and \eqref{eqb04}) in order to reach the equivalent of Eq.\eqref{eqb07} in the case of hypergeometric polynomials having the same behaviour as the $e_n$(See Eq.\eqref{eqf04}). Notice that the statement given below is analogous to the one in Section 2 (See Eq.\eqref{eqb03} and \eqref{eqb04}).

\begin{lemma}
Let us consider the differential equation of the hypergeometric type \eqref{eqa85} with the polynomials $A$ and $B$ given by \eqref{eqf01} and the constant $\lambda$ defined in \eqref{eqa106}. If the condition \eqref{eqf04} is respected, then there is a function $u$ given by
\begin{equation}
\label{eqa127}
u(x) = \int_{(C)} \frac{\sigma(s)}{(s-x)^{n+1}} ds,
\end{equation}
such that the equation \eqref{eqa85} admits particular solutions of the form :
\begin{eqnarray}
\begin{split}
\label{eqa126}
e_n(x) = \gamma_n\frac{A^n(x)}{\sigma(x)} \ \int_{(C)} \frac{\sigma(s)}{(s-x)^{n+1}} ds,\\
\gamma_n=\frac{n!}{2\pi i}\beta_n.
\end{split}
\end{eqnarray}
\end{lemma}

Equipped with the result \eqref{eqa126}, we can show that the $e_n$ polynomials are solutions of the differential equation \eqref{eqa85} and we can also verify the expression of the generating function of the hypergeometric polynomials $e_n$.\\
Let us now give the theorem allowing to determine the generating function of hypergeometric polynomials of the same type as $e_n$ :

\begin{theorem}
Let us consider the case of hypergeometric polynomials $e_n$ as evoked by Lemma 1. One can show that their generating function $E$ is given by :
\begin{equation}
\label{eqf05}
E(x,t)=\frac{\sigma(s)}{\sigma(x)}\Big\vert_{s=\xi(x,t)},
\end{equation}
with the function $\sigma$ given by \eqref{eqzz03} and $s=\xi(x,t)$, solution of $s-x-A(x)t=0$.\\
\end{theorem}
The relation \eqref{eqf05} applied to the $e_n$ polynomials allows to find the generating function \eqref{eqa72}.\\

\proof{To prove the above theorem, we use the expression of the polynomials $e_n$ given by \eqref{eqa126} and the general definition of the generating functions \eqref{eqb06}, which leads to
\begin{equation}
\label{eqg07}
E(x,t) = \frac{1}{2\pi i \sigma(x)} \int_{(C)} \frac{\sigma(s)}{s-x} \ \sum_{n=0}^{+\infty} \left(\frac{A(x) t}{s-x}\right)^n ds.
\end{equation}
Using the fact that the serie in \eqref{eqg07} is geometric, we obtain :
\begin{equation}
\label{eqg07b}
E(x,t) = \frac{1}{2\pi i \sigma(x)} \int_{(C)} \frac{\sigma(s)}{s-x-A(x)t} ds.
\end{equation}
The value of the integral \eqref{eqg07b} can be obtained by using the residue theorem. Indeed, the pole is the solution $s=\xi(x,t)$ of the equation $s-x-A(x)t=0$. The integral representation \eqref{eqa126} is often known, in the literature, as "the Schlaefli integral"\cite{Art03-11}. It coincides with the polynomial representation if the contour $(C)$ surrounds the pole, and must be such that the function $\sigma$ is analytical everywhere on and in $(C)$(See \cite{Art03-11}). This is the case for the pole $s=\xi(x,t)$. So, by virtue of the residue theorem, we obtain the formula \eqref{eqf05}.}

%%%%%%%%%%%%%%%%%%%%%%%%%%%%%%%%%%%%%%%%%%%%%%%%%%%%%%%%%%%%%%%%%%%%%%%%%%%%%%%%%%%%%%%%%%%%%%%%%%%%%%%%

\section{Complete Solution of the Hypergeometric Differential Equation \eqref{eqa75}}

\subsection{General Solution of the Equation \eqref{eqa75}}
\

The differential equation \eqref{eqa75} is of second order. It must have two linearly independent solutions. We stated in the previous section that the family of polynomials $e_n$ is one of its solutions. In this section, we will show by other means that it is so and that there is a second linearly independent solution to that differential equation. \\

Consider the differential equation
\[x y'' + \left(x-n\right) y' - n y = 0.\]
Let $y(x)=e^{-x} \ u(x)$; the equation above becomes
\begin{equation}
\label{eqf06}
x u'' - \left(x + n\right) u' =0.
\end{equation}
By integrating \eqref{eqf06}, we find
\[u(x) = C_1 + C_2 \int x^ n \ e^x dx, \]
and
\begin{equation}
\label{eqf07}
y(x) = C_1 \ e^{-x} + C_2 \ e_n(x).
\end{equation}
This result proves that the family of polynomials $e_n$ is a solution of the differential equation \eqref{eqa75}. In addition, the second solution is the function $y=e^{-x}$. We can see the analogy with the results obtained before (See \eqref{eqz06}).

\subsection{The $e_n$ Polynomials and the Associated Laguerre Polynomials "of a Particular Type"}
\

The associated Laguerre polynomials and the $e_n$ polynomials are hypergeometric polynomials. A careful observation shows similarities between their (second-order) differential equations. We are led to think that there is a relationship between these two families. To achieve this relation, we start by comparing the differential equation of the associated Laguerre polynomials(See \cite{Art03-05, Art03-11} and from \cite{Art03-39} to \cite{Art03-44})
\begin{equation}
\label{eqy01}
x\frac{d^2}{dx^2}L_n^{(\alpha)}(x) + \left(\alpha + 1 - x\right)\frac{d}{dx}L_n^{(\alpha)}(x)+n \ L_n^{(\alpha)}(x)=0
\end{equation}
and Eq.\eqref{eqa75}. We notice that two changes must take place on Eq.\eqref{eqy01}.\\
First, we have to change the variable $x$ to $-x$, which leads to the differential equation below :
\begin{equation}
\label{eqy02}
x\frac{d^2}{dx^2}L_n^{(\alpha)}(-x) + \left(\alpha + 1 + x\right)\frac{d}{dx}L_n^{(\alpha)}(-x)-n \ L_n^{(\alpha)}(-x)=0.
\end{equation}
The second change follows from the comparison of the differential equations \eqref{eqa75} and \eqref{eqy02}. From this, we find that the order of the associated Laguerre polynomials must be $\alpha=-n-1 $. Thus, we conclude that the two families of polynomials must be proportional. We write :
\begin{equation}
\label{eqy03}
e_n(x)=\lambda \ L_n^{(-n-1)}(-x),
\end{equation}
with $\lambda$ a non-zero constant.\\ 
This result is not surprising because the two families being hypergeometric, they can be deduced from hypergeometric functions and confluent hypergeometric functions which are general solutions of the hypergeometric differential equation \eqref{eqa85}. Indeed, it is possible to express a number of elementary and special functions in terms of hypergeometric functions and confluent hypergeometric functions \cite{Art03-05, Art03-09, Art03-11, Art03-12, Art03-13, Art03-08, Art03-20}.\\

To determine the proportionality constant, we write the Rodrigues formula associated with 
$L_n^{(-n-1)}(-x)$(See \cite{Art03-05, Art03-09, Art03-11, Art03-13, Art03-35}) :
\begin{equation}
\label{eqy04}
L_n^{(-n-1)}(-x)=\frac{e^{-x}}{n!} \ x^{n+1} \ \frac{d^n}{dx^n}\left(x^{-1} \ e^x\right)
\end{equation}
that we compare with the $e_n$'s Rodrigues formula(See \eqref{eqa123}), we find :

\begin{equation}
\label{eqf17}
e_n(x) = n! \ L^{(-n-1)}_n (-x).
\end{equation}
The relation \eqref{eqf17} between these two families of polynomials can be verified by comparing the generating function of the associated Laguerre polynomials, obtained by H. J. Weber\cite{Art03-16} and that of the $e_n$ polynomials given in \eqref{eqa72}.\\
It is easy to verify this, knowing the expressions of the Laguerre polynomials of order $\alpha$, that from \eqref{eqf17}, we get $e_0(x)=1, e_1(x)=x-1, e_2(x)=x^2 -2x + 2, \ldots $

%%%%%%%%%%%%%%%%%%%%%%%%%%%%%%%%%%%%%%%%%%%%%%%%%%%%%%%%%%%%%%%%%%%%%%%%%%%%%%%%%%%%%%%%%%%%%%

\section{Generalization of the E Problem}
\

The "E problem", as stated by Theorem 1, can be generalized by a homothety on the argument of the exponential.\\
We define the $E_n^{(m)}(x)$ as the sequence of functions defined by the integral :
\begin{equation}
\label{eqh01}
E_n^{(m)}(x)= \int x^n \ e^{mx} dx,
\end{equation}
with $n, m \in \mathbb{N} \times \mathbb{R}^\star.$\\
We have that 
\begin{equation}
\label{eqh02}
E_n^{(m)}(x)=e_n^{(m)}(x) \ e^{mx} + C,
\end{equation}
where $e_n^{(m)}$ are polynomials of degrees $n$ and order $m$($m \in \mathbb{R}^\star$); $C$ is an integration constant.\\
Obviously, the generalized polynomials $e_n^{(m)}$ and the polynomials $e_n$ are linked. We have :
\begin{equation}
\label{eqk01}
e_n=e_n^{(1)}; \ \forall \ n \in \mathbb{N}.
\end{equation}

Always inspired by the study of the family of the $e_n$ polynomials, we can prove that the polynomials $e_n^{(m)}$ are solutions of the following first-order differential equation :
\begin{equation}
\label{eqh03}
\frac{d}{dx}e_n^{(m)} + m e_n^{(m)} = x^n.
\end{equation}
An explicit expression of these polynomials is found using the Frobenius method:
\begin{equation}
\label{eqh05}
e_n^{(m)} (x) =  m^n x^n + \sum_{l=0}^{n-1} \ \left(-1\right)^{l+n} \ m^l \ \frac{n!}{l!} \ x^l.
\end{equation}

The $e_n$ polynomials, the $e_n^{(m)}$ are hypergeometric. They are the polynomial solutions of the following hypergeometric differential equation :
\begin{equation}
\label{eqh04}
x \frac{d^2}{dx^2}e_n^{(m)} + \left(mx-n\right)\frac{d}{dx}e_n^{(m)} - mnP_n^{(m)}=0.
\end{equation}

Using the NU method, we can find the adapted Rodrigues formula for these polynomials. The fonction $\rho$, solution of Eq.\eqref{eqi01} is given by
\begin{equation}
\label{eqh06}
\rho(x)=x^{-n-1} \ e^{mx}.
\end{equation}
Then the Rodrigues formula is
\begin{equation}
\label{eqh07}
e_n^{(m)}(x)=x^{n+1} \ e^{-mx} \frac{d^n}{dx^n} \left( x^{-1} \ e^{mx} \right).
\end{equation}
Some explicit expressions are : $e_0^{(m)}(x)=1, \ e_1^{(m)}(x)=mx-1, 
\ e_2^{(m)}(x)=m^2x^2 -2mx+2, \ldots$\\

Finally, he generating function $E^{(m)}(x,t)$ of $e_n^{(m)}$, which can be determined by the two procedures(establishment of a differential equation involving it or using the relation \eqref{eqf05}), is :
\begin{equation}
\label{eqh08}
E^{(m)}(x,t)=\frac{e^{mxt}}{1+t}.
\end{equation}

%%%%%%%%%%%%%%%%%%%%%%%%%%%%%%%%%%%%%%%%%%%%%%%%%%%%%%%%%%%%%%%%%%%%%%%%%%%%%%%%%%%%%%%%%%%%%%
%%%%%%%%%%%%%%%%%%%%%%%%%%%%%%%%%%%%%%%%%%%%%%%%%%%%%%%%%%%%%%%%%%%%%%%%%%%%%%%%%%%%%%%%%%%%%%

\section*{Conclusions}
\

This work was aimed at constructing families of polynomials to solve the SCE problems. After investigation, we found relationships between the three problems. Thus, the search for these polynomial sequences can simply be summarized in the determination of only one of them. The sequence $\lbrace e_n; n \in \mathbb{N} \rbrace$ has some interesting properties.\\
First, they are hypergeometric polynomials. As a result, the NU method can be used to demonstrate or verify some results : Rodrigues formula, generating function, ... As for the generating function, the formula \eqref{eqb07} of Nikiforov-Uvarov can not be used without care. Thus, we have established a new result (lemma 1) taking into account the peculiarity of $e_n$ polynomials. This property is the singular presence of the integer $n$ in the expression of the polynomial $B$ of the differential equation (See Eq.\eqref{eqa75}). The corollary of this property is that the function $\sigma(x)\equiv\rho_n(x)$(See Eq.\eqref{eqzz02} and \eqref{eqzz03}), for the $e_n$, is independent of the integer $n$. This result does not occur for classical hypergeometric polynomials (See Eq.\eqref{eqf03} and \eqref{eqf04}).\\
Second, the $e_n$ polynomials can be deduced from hypergeometric functions. Indeed, the differential equation \eqref{eqa85} is a confluent hypergeometric equation \cite{Art03-11, Art03-12, Art03-13}. This result leads to a relation between the $e_n$ polynomials and the associated Laguerre polynomials with negative orders.\\

The use of the Frobenius method for the determination of the three families of polynomials, related to SCE problems, made it possible to find the results of G. Dattoli \textit{et al} \cite{Art03-36, Art03-37, Art03-38} and H. M. Srivastava \textit{et al} \cite{Art03-46}. They've obtained expressions of the integrals \eqref{eqa04},\eqref{eqa05},\eqref{eqa06} in terms of polynomial series from a negative derivative operator and the truncated exponential. Our approach is different in the sense that we began by establishing and demonstrating Theorem 1 expressing the integrals \eqref{eqa04}, \eqref{eqa05} and \eqref{eqa06}. The rest of the work focused on the link with the NU method, the determination of the Rodrigues formula and the study of the particular properties of the hypergeometric polynomials of the same type as the $e_n$. Indeed, the generating functions of such families can not be deduced, directly, from the result \eqref{eqb07} of the NU method. This last fact led to the development of Theorem 2 adapted to the situation.\\

At the end of our investigation, we obtained four procedures for determining the $e_n$ polynomials : the first is their expressions in the form of polynomial series given by \eqref{eqa64}. This required the use of the Frobenius method in solving the differential equation \eqref{eqa03}. The second is the Rodrigues formula \eqref{eqa123}, obtained using the NU method. The third is the expression of the generating function \eqref{eqa72}. This is obtained by two methods : the resolution of the differential equation of which it is a particular solution and the use of the lemma 1(the adapted NU method in this case). This theorem is inspired by the NU method. The last one is the expression of the $e_n$ as a function of the associated Laguerre polynomials with negative order (See Eq.\eqref{eqf17}). The present study shows how some mathematical problems can lead to results connected to other ones obtained by various methods and procedures.\\

This study showed that the generating function has to be computed with care when one of the polynomial ($A$ or $B$) appearing in the hypergeometric equation involves an integer $n$ such that the function $\rho_n$(See \eqref{eqzz02}) does not depend on this integer. We explicitly gave a recipe which works in the case of the SCE problems. Everything was in control since the results we presented could be retrieved by other methods. The case we dealt with involved a $B$ function which is linear in $n$. A natural generalization going above the SCE problem will treat non trivial dependences of the funtions $A$ and $B$ on a natural integer $n$. The path followed here already brings fruitfull results. We are working on it and hope to present the corresponding results soon.

%%%%%%%%%%%%%%%%%%%%%%%%%%%%%%%%%%%%%%%%%%%%%%%%%%%%%%%%%%%%%%%%%%%%%%%%%%%%%%%%%%%%
%%%%%%%%%%%%%%%%%%%%%%%%%%%%%%%%%%%%%%%%%%%%%%%%%%%%%%%%%%%%%%%%%%%%%%%%%%%%%%%%%%%%

\section*{Acknolegments}
\

We would like to express our sincere thanks to Professors Butsana-Bu-Nianga and Musongela Lubo of the Department of Physics of the Faculty of Sciences of the University of Kinshasa. The first one for initiating us to the SCE problem and leading our first investigations. The second one for his keen scrutinizing eye on this work. 

%%%%%%%%%%%%%%%%%%%%%%%%%%%%%%%%%%%%%%%%%%%%%%%%%%%%%%%%%%%%%%%%%%%%%%%%%%%%%%%%%%%%%%
%%%%%%%%%%%%%%%%%%%%%%%%%%%%%%%%%%%%%%%%%%%%%%%%%%%%%%%%%%%%%%%%%%%%%%%%%%%%%%%%%%%%%%

\end{document}